\documentclass{amsart}
\usepackage{latexsym,amsxtra,amscd,ifthen}
\usepackage{amsfonts}
\usepackage{verbatim}
\usepackage{amsmath}
\usepackage{amsthm}
\usepackage{amssymb}

\numberwithin{equation}{section}

\theoremstyle{plain}
\newtheorem{theorem}{Theorem}[section]

\newtheorem{prop}[theorem]{Proposition}

\newtheorem{cor}[theorem]{Corollary}

\theoremstyle{definition}
\newtheorem{defi}[theorem]{Definition}
\newtheorem{definition}[theorem]{Definition}

\newtheorem{rema}[theorem]{Remark}
\newcommand{\s}{\ensuremath{\mathbb S}}
\newcommand{\Z}{\ensuremath{\mathbb Z}}

\def\char{\mbox{char\,}}

\def\dim{\mbox{dim\,}}
\def\id{\mbox{id\,}}

\def\im{\mbox{im\,}}

\def\ker{\mbox{ker\,}}

\def\Palg{{{\mathcal P}\mbox{-Alg\,}}}

\def\tor{\mbox{Tor\,}}

\def\ass{\: {\mathcal A}ss \:}
\def\com{\: {\mathcal C}om \:}
\def\gcom{\: {\mathcal G}en{\mathcal C}om \:}

\def\p{ {\mathcal P} }

\begin{document}

\title{On Kurosh problem in varieties of algebras%
}\thanks{Partially supported by the Dynastia
foundation, by  the grant of President of Russian Federation
MD-2007.288.1,  by the Program of Support of Scientific Schools grant 1983.2008.1, and
by the Russian Foundation of Basic Research grant 08-01-00297-a}%

\author{Dmitri Piontkovski}

      \address{Department of High Mathematics for Economics,
Myasnitskaya str. 20, State University `Higher School of Economics', Moscow 101990, Russia
% Central Institute of Economics and Mathematics\\
%                      Nakhimovsky prosp. 47, Moscow 117418,  Russia
}

\date{December 1, 2008}

\email{piont@mccme.ru}



%\date{\today}

 \begin{abstract}
 We consider a couple  of versions of  classical Kurosh problem
 (whether there is an infinite-dimensional algebraic algebra?)
 for varieties of linear  multioperator algebras over a field.
 We show that, given an arbitrary signature, there is a variety
 of algebras of this signature such that  the free algebra of the
 variety contains multilinear elements of arbitrary large  degree,
  while the clone of every such element satisfies some nontrivial
  identity. If, in addition, the number of binary operations is
  at least 2, then one can guarantee that  each such  clone is finitely-dimensional.

Our approach is the following: we translate the problem to the
language of operads and then apply usual homological constructions
in order to adopt  Golod's solution of the original Kurosh
problem.

The paper is expository, so that some proofs are omitted. At the
same time, the  general relations of operads, algebras, and
varieties are widely discussed.
 \end{abstract}

\maketitle

%\section{A dictionary of universal algebra: \\
%operads vs identities}

\section{Introduction}

\subsection{Kurosh problem and multioperator algebras}

In his prognosis of the future of general algebra, A.~G.~Kurosh
predicted in 1970 that the main interests of general algebra would
moved, in the nearest decades, to the neutral territory between
universal and general algebra~\cite[p.~9]{kurlec}. Following this
tendency, in  this note we try to extend some classical algebraic
ideas of the time of Kurosh to a territory very close to universal
algebra.

We start with a classical question of Kurosh: whether there can be
a finitely-generated infinite-dimensional algebraic algebra? The
first solution was found by  E.~S.~Golod~\cite{golod} in 1964, who
constructed a number of analogoues examples for different
algebraic systems (such as associative algebras, Lie algebras and
$p$--groups). Here we will discuss the Kurosh problem in the
setting of linear universal algebra, for algebraic systems in
place of algebras in the formulating. To construct suitable
examples, we extend the original Golod technique, including the
Golod--Shafarevich theorem~\cite{gs}.

We consider the varieties of multioperator linear
$\Omega$--algebras over a field $k$ (the term
from~\cite[\S~13]{kurlec}). Every multilinear element of the
countable generated free algebra $F$ of a variety  $W$ may be
identified with a family of multilinear operators (or operations)
acting on the algebras of $W$. The linear combinations of finite
compositions of a multilinear operator with itself form a linear
subspace of $F$ called {\it clone} of the operator. We will see
%in Section~\ref{sec:oper_n_alg}
that the notion of clone of a single operator is analogous to the
concept of  one-generator subalgebra of an associative algebra.
Therefore, one can consider two natural versions of the above
Kurosh problem. The strong version (we will refer it as Burnside
problem) is the following: {\it does there exist a variety $W$ (of
a given finite signature  $\Omega$) such that it admits nonzero
$n$-linear operations for all arbitrary large $n$ while every
clone of a single operation is finite-dimensional}? A weaker
version of the problem is the following: {\it does there exist a
variety $W$, again with infinite-dimensional space of multi-linear
operations, such that every such an operation satisfy nontrivial
``identity'', that is, some nontrivial linear combination of
 composition of the operation with itself is zero}?

 One can consider also relative versions of the above problems, that is, wuth
 an additional restriction that  the variety $W$
must be a subvariety of a given variety ${\mathcal V}$. In this
paper, we focus on the non-relative versions, but our method
(based on Theorem~\ref{th:GS_oper}) gives also some examples to
the relative case.

\subsection{The language of operads}

 In linear universal algebra (when the algebras are modules over a ring or a field and the
 operations are multi-linear homomorphisms), there are two main
languages for definitions and theorems which are used for different purposes.

First one is the classical language of {\it varieties} and {\it
identities}: it is the most popular one in ring theory and in
general algebra as well.  The theory of algebras with polynomial
identities is usually described in this language. The works of
Kurosh are also  written in it.

Another language is based on the concept of  {\it operad}. It is
used mostly in algebraic topology and mathematical physics. This
language seems more appropriate for discussions on the homological
properties of algebraic systems, including Massey operations and
Kontsevich formality.

Our first step here is to give a brief dictionary of these two
languages  in Section~\ref{sec:oper_n_var}. In this dictionary, we
define operads and related notions in terms of varieties; dually,
we define the varieties, identities etc in operadic terms.
Therefore, it is assumed that the reader do understand at least
one of these two languages. We hope that this brief dictionary
will help the both kind of universal algebraists to understand
each other,
--- or, at least, to recognize the well-known objects in foreign-language
descriptions. We use this dictionary in order to get a more
natural formulation of our versions of universal Kurosh problem.

For the reader convenience, let us give a rough and brief
translation table in a phrase-book style.

\medskip
\centerline{\bf A phrase-book}
\smallskip
 \centerline{
 \begin{tabular}{rcl}
  variety & --- & operad \\
  subvariety & --- & quotient operad \\
  clone & --- & suboperad \\
  signature &
  % of a variety
  --- & set of generators \\ % of an operad
  identities  & --- & relations \\
  free algebra & --- & free algebra \\
  codimension series & ---  & generating function \\
    T-space & --- & right ideal \\
  T-ideal & --- & ideal \\
   Specht properties & --- & Noether properties\\
  % Kemer theorem & --- & Hilbert Basis theorem \\
   \end{tabular}
  }

  \medskip

Now, our two versions of Kurosh problem are translated as follows.
Burnside problem (or strong Kurosh problem)  for operads: {\it
does there exist an infinite finitely generated operad $P$ such
that every its one-generated suboperad is finite}? Kurosh problem
(in a weak version) for operads: {\it does there exist an infinite
finitely generated operad
 $P$ such that every its one-generated suboperad is not free}?
 The relative version sounds as follows: {\it given an operad ${\mathcal S}$,
 can one choose the operad $P$  above to be a quotient of ${\mathcal S}$}?

In this paper, we give partial answers to both these problems.
That is, we show that such varieties and operads do exists (in the
stronger version, with a restriction on the signature). The
relative versions are in general open.

Let us give a ``bi-lingual'' formulation (some details are given
in Section~\ref{sec:gs} below).

\begin{theorem}[Corollaries~\ref{cor:weak_kur},~\ref{cor:strong_kur}]
1) Let $\Omega$ be a finite signature. Then there exists a variety
of algebras of signature $\Omega$ such that there are nonzero
multilinear operations of arbitrary high order but the clone of
each operation in a free algebra satisfies a non-trivial identity.
If $\Omega(2)$ has at least 2 elements, then one  can assert, in
addition, that the clone of each operation in a free algebra is
finite-dimensional.

2) Let $X = X(2), X(3), \dots $ be an $\s$-module, that is,  a
sequence of representation of symmetric groups $S_2, S_3, \dots$
Then there exists an operad $\p$ generated by $X$ such that each
 one-generated suboperad in it is not absolutely  free. If $\dim X_2
\ge 3$, then we can assert, in addition, that every its
one-generated suboperad is finite.
\end{theorem}

\subsection{Organization of the paper}

{\it This paper is expository; most proofs are omitted and will be
published in  the subsequent paper}~\cite{to_appear}.

The dictionary mentioned before is
given  in  Section~\ref{sec1}. % we give the dictionary.
%We conclude it with the translation of the  Birkhoff HSP theorem
%to the operadic language.
In  Section~\ref{sec2}, we discuss the
(well-known) analogy between  operads and  graded associative algebras. This leads to an analogy
of linear universal algebra and graded ring theory. So, we describe  the notions
of ideals, modules, generators and other algebraic  terms in operad theory.

This analogy allows to develop a version of classical homological
algebra, including free resolutions and derived functor, in the
category  of modules of given operads. We use it in order  to
transfer the Golod--Shafarevich theorem to operads, see
Section~\ref{sec:gs}. This leads, by the Golod method, to a
construction of infinite operad with finite one-generated
sub-operads. This gives the solutions of weak and strong Kurosh
problems for varieties (operads) for  algebras of almost arbitrary
signature,   see Corollaries~\ref{cor:weak_kur}
and~\ref{cor:strong_kur}.

%As an other  application of our
%approach, we solve a version of Burnside-type problem for %operads
%varieties of linear algebras (or, equivalently, for operads) by
%using an operadic version of the famous Golod-Shafarevich theorem.

\subsection{Assumptions}

We consider varieties of multioperator linear algebras over a
field $k$. To avoid technical details, we assume that  the
signature $\Omega$ of a variety is always locally finite and does
not contain constants and unary operations, that is, $\Omega =
\Omega_2 \cup \Omega_3 \cup \dots $ there the subsets $\Omega_n$
of $n$-ary operations (i.e., $n$-linear operators) are finite. To
simplify the notation, we assume that the identical operator
(which does not belong to $\Omega$) is also applicable to  any
algebra of the variety. In Section~\ref{sec:oper_n_var} we assume,
unless otherwise is stated, that the ground field $k$ has zero
characteristic.
%Also, we assume that the ground
%field $k$ has zero characteristics.

\subsection*{Acknowledgement}

I am grateful to %Viktor Kac,
Alexei Kanel--Belov and Louis Rowen
%,  and
for discussions on the definitions of varieties and operads, and
to Anton Khoroshkin for useful remarks.
  %here at IHES.
  My special thanks to Evgenij Solomonovich Golod, who had
introduced me, years ago, to the world of algebraic operations.

I am grateful to Bar--Ilan University, where the paper had been
partially wrote, for hospitality and friendly atmosphere. I am
grateful to  the Gelbart Institute and the Bar-Ilan Mathematics
Department for funding my visit to Bar--Ilan.

\section{Operads vs varieties: a dictionary}

\label{sec1}
\label{sec:oper_n_var}

\subsection{A definition of operad}

Let $W$ be a variety of $k$--linear algebras (without constants,
with identity and without other unary operations)  of some
signature $\Omega$. Consider the free algebra $F^W (X)$ on a
countable set of indeterminates $X= \{ x_1, x_2, \dots \}$. Let
$\p (n) \subset F$ be the subspace consisting of all multilinear
generalized homogeneous polynomials on the variables $x_1, \dots,
x_n $.

\begin{defi}
Given such a variety $W$,  the sequence $\p_W = \p := \{ \p (1),
\p (2), \dots \}$ of the vector subspaces of $F^W (X)$ is called
an {\it operad}\footnote{More precisely, symmetric connected
$k$--linear operad with identity.}. The signature $\Omega$ is
called a {\it generation set} of the operad $\p_W$.
\end{defi}

The $n$-th component  $\p (n)$ may be identified with the set of
all derived $n$-linear operations on the algebras of $W$; in
particular,  $\p (n)$ carries a natural structure of a
representation of the symmetric group $S_n$. Such a sequence $Q =
\{ Q(n) \}_{n \in \Z}$ of representations $Q(n)$ of the symmetric
groups $S_n$ is called an {$\s$--module}, so that an operad
carries a structure of $\s$-module.  Also, the compositions of
operations (that is, a substitution of an argument $x_i$ by a
result of another polylinear operation, with a subsequent
re-numerating the variables) gives natural equivariant maps of
$S_*$-modules $\circ_i :  \p(n) \otimes \p(m) \to \p(n+m-1)$. Note
that the axiomatization of these operations gives an abstract
definition of operads, see~\cite{oper} for the discussion on
different definitions.

A morphism of operads $f: \p \to \p'$ is a sequence of maps $f(n): \p
(n) \to \p' (n)$ of $S_n$-modules compatible with  the compositions.
For example, any inclusion $W \subset W'$ of varieties of the same
signature gives a surjective operadic  map $\p_W \to \p_{W'}$.

\subsection{A definition of variety}

Let $\p$ be an operad (with $\p(0)$ equal to zero and
one-dimensional $\p(1)$ spanned by the identity element) with some
discrete generating set $\Omega$. Recall that an {\it algebra}
over $\p$ (=$\p$-algebra) is a (non-graded) right $\p$-module,
that is, a vector space $V$ with $\p$-action $ \p(n): V^{\otimes
n} \to V$ compatible with compositions of operations and the
$k[S_*]$-module structures on $\p$. The class (or an additive
category) of all algebras over $\p$ is denoted by $\Palg$.

\begin{defi}
Let  $\p$ be an operad with generating set  $\Omega$.
 Then the category $\Palg$
%class $W=W_{\p}$ of all algebras over $\p$
is called a {\it variety} of algebras of signature $\Omega$.
\end{defi}

{\it From now, let us fix a variety $W$ and a corresponding operad
$\p$, so that $W= W_\p$ and $\p = \p_W$. We fix also a minimal
generating set (or signature) $\Omega$}. The implications $(W= W_\p)
\Longleftrightarrow (\p = \p_W)$ will be discussed later in Proposition~\ref{linearization}.

\subsection{(Co)dimensions}

An $n$-th codimension of a variety $W$ is just the dimension of the
respective operad component: $c_n(W) = \dim_k \p_W(n)$. The {\it
codimension series} of the variety $W$, or the {\it generating
function} of the operad $\p$, is  a formal power series
$$
\p (z) := \sum_{n \ge 1} \frac{c_n(W)}{n !} z^n = \sum_{n \ge 1}
\frac{\dim_k \p_W (n)}{n !} z^n .
$$
An analogous  generation function $Q(z) = \sum_{n \ge 1}
\frac{\dim_k Q (n)}{n !} z^n$ is defined by for every $\s$--module
$Q$ with $\dim_k Q(n) < \infty$ for all $n \ge 1$.
 % where $z$ is a formal variable.
For a more general versions of this generating function (which
involve, in particular, the characters of the representations
$\p(n)$ of groups $S_n$) we refer the reader to the well-known
paper of Ginzburg and Kapranov~\cite{gk}.

If the set $\Omega $ is finite, then the series $\p (z)$ define an analytical function
in a neighborhood of zero.
For example, the operad $\ass$ of associative algebras has generating function
$\ass (z) = \frac{z}{1-z}$. For every proper quotient operad $\p$ of $\ass$, we have
 $\lim\limits_{n\to \infty} \frac{\ln
\dim \p (n)}{n} = \ln c(\p)$, where $c(\p) \in \Z$ (due to
Giambruno and Zaitsev, see~\cite{giza} and references therein): in
particular, the function $\p(z)$ is in this case analytical in the
whole complex plane.

%\section{Free operads and free varieties}

% \section{Identities vs relations}

\subsection{Free variety}

Recall that a sequence $M = \{ M(1), M(2),  \dots \}$, where each
$M(i)$ is a $k[S_i]$-module, is called an $\s$-module. Given a
sequence $\Omega = \{ \Omega (1), \Omega (2), \dots \}  $ of
discrete sets, we naturally define a free $\s$-module $\s \Omega =
\{ k[S_1] \Omega (1), k[S_2] \Omega (2), \dots \}  $. In
particular, every operad is an $\s$-module, and every subset of an
$\s$-module generate an $\s$-submodule. If $\Omega$ is a minimal
generating set of an operad $\p$, then the $\s$-submodule  $\s
\Omega$ is also called a generating space of $\p$.

\begin{defi}[of free variety]
The variety $W$ of signature $\Omega$ is called free, if
 the generating set $\Omega$  minimally generates a free $\s$-submodule in the
 operad $\p = \p_W $ and the operad $\p$ itself is  free with generating $\s$-submodule $\s
\Omega$.

A free algebra in a free variety is called absolutely free (of given
signature).
\end{defi}

We will call the operad of a free variety {\it absolutely free}.

\subsection{Free operad}

Let $\p$  be an operad generated by a subset $\Omega$. The operad
$\p$ is called free (on the generating set $\Omega$), if the
T-ideal $T$ of identities of the variety $W = W_\p$
%(or,  equivalently, the ideal $I$ of relations of $\p$)
consists of the  linear combinations of generators, that is, by
elements of the $\s$-submodule $X$ of $\p$ generated by  $\Omega$.
Since the free operad $\p$  is uniquely determined by the
$\s$-submodule $X$, it is denoted by $\Gamma ( X) $ (notation
from~\cite{oper}).

For example, any absolutely free operad is free (since it has no relations).
On the other hand,
the operad $\gcom$ of general (nonassociative) commutative algebras
is free with the multiplication $\mu$ as a generator ($\Omega = \{ \mu \}$),
but is not absolutely free because of the identity $[x_1, x_2]:= \mu(x_1, x_2) - \mu \circ \tau (x_1, x_2) \in X$,
where $\tau$ is the generator of the group $S_2$.

\subsection{Relations of operads}

For every two $\s$-submodules $A$ and $B$ of an operad $\p$, one
can define a new  $\s$-submodule $A\circ B \subset \p$ generated
by all compositions $a(b_1, \dots, b_n)$, where $a \in A \cap
\p(n)$ and $b_i \in B$. An $\s$-submodule $I \subset \p$ is called
a left (respectively, right, two-sided) {\it ideal,}  if $I = \p
\circ I$ (resp., $I = I \circ \p$, $I = \p \circ I \circ \p$). The
generating sets of ideals are defined in the obvious way.

It follows that the two-sided ideals are exactly the kernels of
operadic morphisms. If an operad $\p$ is represented as a quotient
(=`image of a surjective morphism') of a free operad $\p'$ by a
two-sided ideal $I$, the elements of $I$ are called the {\it
relations} of  $\p$. Given a generating set $\Omega$ of $\p'$, all
the relations becomes the identities of the variety $W_{\p'}$ in
this signature.

For example, the operad $\com$ of commutative associative
algebras, as a quotient of the free operad $\gcom$  described
above, has the associativity relation $Ass(x_1, x_2, x_3):= (x_1
\cdot x_2) \cdot x_3 -x_1 \cdot (x_2 \cdot x_3)$ (where $a\cdot
b:= \mu(a,b) $), and all other relations belong to the two-sided
ideal $I$ in $\gcom$ generated by this relation $Ass(x_1, x_2,
x_3)$.

\subsection{Identities of varieties}

Let $F$ be a free algebra  over an operad $\p$ . Then one can
consider every  element of $F$ as an operation on other algebras
of the variety, where the generators of $F$ are replaced by
elements of another algebra. So, for every  $p, q \in F$ one can
define a composition $p\circ_i q \in F$ (it is equal to $p$ if $p$
does not really depend on the variable $x_i$). Note that, in
contrast to the operad composition, there is no re-numerating of
variables after c composition, e.~g.,  $p(x_1,x_2) \circ_1 x_2 =
 p(x_2,x_2)$ etc.  Then for every subset $C
\subset F$ one can define two composition subset $F \circ C = \{ f
\circ_i c \}$ and $C \circ F = \{ c \circ_i f \}$ (where $f,c,i$
runs through $F, C$ and $\Z_{\ge 1}$, respectively). A linear
subspace $C \subset F$ is called an ideal (respectively, T-space,
subalgebra) if it is closed under composition of the first type
(respectively, second type, both types).

Let $W'$ be a free variety with signature $\Omega$, and let $\p' =
\p_{W'}$ be the corresponding operad.  The ideal $I$ of relations of
the operad $\p$ in $\p'$ generates a T-ideal $T$ in the absolutely
free algebra $F = F^{W'}(X)$ with a countable generating set $X$.
The elements of $T$ are called {\it identities} of the variety $W$.

The standard linearization process gives a procedure to establish
the following

\begin{prop}
\label{linearization} Let $F$ be a free algebra with countable
generating set  of a variety $W$ over a field $k$ of zero
characteristic. Then every T-space or a T-ideal $Y$ in $F$ is
generated by the subset $Y \cap \p_W$ of multilinear elements.
Moreover, this subset $Y \cap \p_W$ form an ideal (right- or
two-sided, respectively) in the operad $\p_W$.

In particular, it follows that the T-ideal $T$ of identities of an
arbitrary variety $W$ is generated by the relations of the operad
$\p := \p_W$, hence $W = W_\p$.
\end{prop}

The proof of this proposition (essentially, it is a description of
the linearization process mentioned above) is essentially the same
for all types of linear algebras; see, e.~g.,~\cite{bero}
or~\cite{giza}
 for the case of associative PI algebras.  For example, the linearization of an identity $x_1^2 $
 in the free algebra of $W_{\gcom}$ leads to the identity $x_1 x_2$
 (because $x_1 x_2 = \frac{1}{2} ((x_1+x_2)^2 - x_1^2 - x_2^2)$),
 that is, the identity $x_1^2$  defines the variety of algebras with zero
 multiplication, that is, the category of vector spaces.

 Note that the linearization is essentially depend on the assumption $\char k = 0$.
 If the characteristic of the field is positive,
 then only the implication $(W = W_\p) \Longrightarrow (\p = \p_W)$ is valid,
 but the reverse implication fails.

%The operad $\Com$ of commutative algebras

\section{Operads \& graded algebras:  \\an analogy }

\label{sec:oper_n_alg}

 \label{sec2}

A homological theory of operads is similar to the homological
theory of associative algebras.
 There is a number of homological construction, which are successively
 moved
 from rings to operads: (co)bar constructions, minimal models, (DG-)resolutions and (DG-)modules,
 Koszul duality etc.~\cite{fresse, gk, oper}.
Here we move to operads a part of classical homological algebra,
namely, the theory of torsion functors. By a standard way, we will
construct free resolutions of modules over operads and use them to
define and calculate  the derived functors of an operadic analogue
of tensor product. This will be used later in our version of
Golod--Shafarevich theorem.

  Note that for every operad $\p$ one can define graded right and left modules
 over it: they are $\s$-modules $V$ with the structure of $\p$-algebras (right modules)
  or with the compositions $V(n) \circ_i \p(m) \to V(n+m-1)$  (left modules),
  where in both cases  the structure should be compatible     with the operadic and $\s$-module structures.
 The composition functor
    $- \circ_{\p} L $ (where $L$ is a graded left $\p$-module) from the category $mod-\p$ of graded right
    $\p$-modules to the category $k-mod$ of graded vector spaces over $k$
    is analogous
    to tensor product of modules over a graded algebra.  It has left derived functors $\tor^\p_i (R,L)$
    which are analogous to usual $\tor$s of modules over graded algebras. These operadic
    torsion functors can be calculated using free resolutions (or cofibrant resolutions, in the DG case)
    of the first argument~\cite{fresse}.

    A formal explanation of these ideas can be given by the following
    chain of standard statements.

    \begin{prop}
    Let $\p$  be an operad. % with minimal generating set  $\Omega$
    Then the category $mod-\p$  of all graded  right $\p$-modules
    is abelian
    (hence, it is an abelian subcategory of the category of all
    $\s$-modules).
    \end{prop}

    Let $V$ be an arbitrary  $\s$--module. A composition right $\p$--module $V \circ \p$
    is called {\it free} (and $V$ is called as its minimal $\s$-module of generators).
    As an $\s$-module, it is a composition product, so that its
    generation function is equal to $(V \circ \p) (z) = V(\p(z))$.
    For example, $\p$ itself is a free right module generated by
    the trivial $\s$--module
    $k \, \id $.

    Suppose that  $M$ is a right graded $\p$--module minimally generated by
    an $\s$--module $V'$ isomorphic to $V$. Then there is a
        (unique up to an isomorphism  $V \to V'$)
    surjective map of $\p$-modules $p: V \circ\p \to
    M$ which isomorphically maps $V$ to $V'$. The kernel $\ker p$
    belongs to the `submodule of decomposables' $V\circ_\p \p_+ \subset  V \circ\p
    $, where $\p_+ = \p(2) \oplus \p(3) \oplus \dots$ is the
    maximal ideal of $\p$. Iterating this construction, we get an
    exact sequence of right graded $\p$-modules
    $$
          \dots F_2 \stackrel{d_2}{\longrightarrow}
                     F_1 \stackrel{d_1}{\longrightarrow}
                    F_0 \stackrel{d_0}{\longrightarrow}
                    M \to 0,
    $$
    where all modules $F_i$ are free (that is why we refer the subsequence
        ${\bf F}: \dots F_2 \to
                     F_1 \to
                     F_0$ as {\it free resolution}
    of $M$). In addition, we have $\im d_i \subset F_{i-1} \circ \p_+$
    for all $i\ge 0$. A resolution with the last property will be called {\it minimal}.

    The second part of the following standard proposition  can proved by the same way
    as for graded modules over a connected graded associative algebra
    (where all graded projective modules are free).

    \begin{prop}
      (i) Every right $\p$-module $M$ admits a minimal free
      resolution ${\bf F}$.

      (ii) If ${\bf F}'$ is another free resolution of $M$, then
      the complex ${\bf F}$ of $\p$-modules is a direct summand of ${\bf F}'$.

      (iii) The minimal resolution ${\bf F}$ is unique up to
      isomorphism of complexes of right $\p$-modules.
    \end{prop}

    Let $R$ and $L$ be right and left graded $\p$--modules. Then
    one can define a composition $k$--module
    $R \circ_\p L$. It is a quotient $k$-module of $ R \circ L$
    by the relations induced by the action of $\p$.

    \begin{prop}
    \label{tors_n_resols}
       Let $L$ be a left graded $\p$--module.

       (i) The functor $C_L: X \mapsto X \circ_\p L $ is right exact on $\p$-mod.

       (ii) There exist a  derived functors $L_* C_L(M)$
       whose value $\tor_i^\p (M,L):= L_i C_L(M)$ can be calculated, for each $i > 0$,
       as the $i$-th homology $\s$-module of the
       complex
       $$
       {\bf F} \circ_\p L: \dots F_2  \circ_\p L \to    F_1  \circ_\p L\to
                     F_0  \circ_\p L.
       $$
    \end{prop}

    Note that the above derived functors has been introduced by
    Fresse~\cite[2.2.4]{fresse} in a more general context of
    DG-operad. The second part of the above
    Proposition~\ref{tors_n_resols} follows from~\cite[Proposition~2.2.5]{fresse}.

     For example, the minimal  $\s$-submodules that generate the
     modules $F_i$ from the minimal free resolution can be calculated
     as
     $$F_i/ F_i \circ \p_+ = \tor_i^\p (M,\p/ \p_+).$$

\section{A criterion for infinite operads and Kurosh problem}

\label{sec:gs}

We give here an operadic version of famous Golod--Shafarevich
theorem which gives a criterion for an associative algebra to be
infinite-dimensional. Despite the original version of the proof of
Golod and Shafarevich~\cite{gs} can be almost directly translated
to the language of operads (with Shafarevich complex replaced by
the first step of the construction of a minimal model of the
operad $\p$, see~\cite{oper}), we prefer to adopt another  approach (explained by
Ufnarovski~\cite{ufn}) based on a direct construction of the minimal free
resolution of the trivial module.

\begin{theorem}
\label{th:GS_oper}
 Let $\p$ be an operad minimally generated by an
$\s$-module $X \subset \p$ with a minimal $\s$-module of relations
$R \subset \Gamma (X)$. We assume here that both these
$\s$-modules are locally finite, that is, all their graded
component are of finite dimension.

Suppose that the formal power series
$$ %\frac{1}{1 - \frac{X(z)}{z} +  \frac{R(z)}{z}}
\left( 1 - \frac{X(z)}{z} +  \frac{R(z)}{z} \right)^{-1}
$$
%is infinite and its
 has non-negative coefficients. Then  the
operad $\p$ is infinite.
\end{theorem}

\begin{proof}[Sketch of proof]

Consider the trivial bimodule $I = \p / \p_+$ (where $\p_+ = \p(2)
\oplus \p(3) \oplus \dots$ is the
    maximal ideal of $\p$, as before). For the generators of the
    beginning of its minimal free resolution, we have
    $\tor_0^\p (I,I)\cong I$, $\tor_1^\p (I,I)\cong X$ and   $\tor_2^\p (I,I)\cong
    R$. This means that the beginning of the resolution looks as
    \begin{equation}
    \label{resol_k}
        0\to \Omega^3 \to R\circ \p \stackrel{d_2}{\to} X\circ \p \to \p \to
        I \to 0,
    \end{equation}
    where $\Omega^3$ is the kernel of $d_2$.

    Taking the Euler characteristics of the exact
    sequence~(\ref{resol_k}), one can we get an equality of formal power
    series
    $$
    \Omega^3(z) =  (R\circ \p)(z) - (X\circ \p)(z) + \p(z) -I(z).
    $$
    Since the formal power series $\Omega^3(z)$ has nonnegative
    coefficients, we obtain the following coefficient-wise inequality
    $$
         R(\p(z)) - X(\p(z)) +\p(z) - z \ge 0.
    $$
    Manipulations with formal power series (including the Lagrangian inverse)  complete the proof.
\end{proof}

For operads generated by binary operations, one can simplify the above condition.

\begin{cor}
\label{cor:GS_binary}
Let $\p$  be an operad and let $X$ and $R$
be as above. Suppose that  $\p$ is generated by binary operations
(that is, $X=X(2)$). Suppose that the function
$$
            \phi(z) = 1 - \frac{X(z)}{z} +  \frac{R(z)}{z}
$$
is analytical in a neighborhood of zero (it is always the case if
$X$ is finitely generated) and  has a positive real root $z_0$ in
this neighborhood such that $\phi(z_0)' \ne 0$. Then the operad
$\p$ is infinite.
%In particular, if  $\dim_k \p(2) \ge 3$ and the $\s$-module $R$ of relations is generated by at
%most one element in every degree $n\ge 27$, then $\p$ is infinite.
\end{cor}

\begin{cor}[Weak Kurosh problem for multi-operational algebras]
\label{cor:weak_kur}
\label{Bur:varieties}
 Suppose that the ground field $k$ is countable.
 Let $\Omega$ be an arbitrary nonempty countable signature.
Then there is a variety $W$ of algebras
 of signature $\Omega$ such that the clone of  every polilinear operation in this variety
 satisfies some non-trivial identity
 while there are multi-linear elements of the free algebra $F^W (x_1, x_2,
 \dots)$ of arbitrary high degrees.
\end{cor}

\begin{proof}
It is sufficient to prove that the suboperad in  $\p_W$ generated
by an arbitrary single multi-linear operation is not absolutely
free. Let us enumerate by positive integers all elements
(=operations) of degree (=arity)
 $\ge 2$ in the absolutely free operad ${\mathcal F }$ generated by $\Omega$.
 If the suboperad $P $ in ${\mathcal F }$ generated by such an operation $p_i$
 (where $i\in \Z_+$) is not absolutely free (e.~g., if the $\s$-submodule $X$,
generated by  $p_i$ in ${\mathcal F }$ is not free), then its
image in $\p_W$ is not absolutely free as well.

Suppose that the operad $P$ is absolutely free. Let us denote by
$R_i$ the sum of all multi-linear compositions of $N=N_i$ copies
of $p_i$, where the numbers  $N_i$ a chosen so that the degrees
$t_i$ of the elements  $R_i$ increase: $t_1 < t_2 < \dots$
 Then every element $R_i$  is invariant under the action of the symmetric group $S_{t_i}$.
 Therefore, the generating function  of the $\s$-module $R$ generated by all these
 elements  $R_i$ is $R(z)
\le \sum_{n\ge 1} z^{t_n}/{t_n!} \le \sum_{n\ge t_1} z^n/{n!}$
(the last coefficient-wise inequality follows from the
inequalities $\dim R (n) \le 1$ for all $n \ge  1$). If the
numbers  $t_i$ are chosen so that they are sufficiently large
($0<< t_1 << t_2 << \dots$), then Theorem~\ref{th:GS_oper} and the
above estimate imply that the generating function of the operad
$\p$ generated by $\Omega$ with the $\s$-module of relations $R$is
infinite.
\end{proof}

The next claim gives a stronger version of Kurosh problem, that
is, the Burnside problem.

\begin{cor}[Burnside problem for multi-operational algebras]
\label{cor:strong_kur}
\label{cor: bur_strong}
 Suppose that the ground field $k$ is countable.

Let  $X= X(2) \cup X(3) \cup \dots$ be an $\s$-module such that
$\dim X(2) \ge 3 $.
% Suppose that the dimension of the representation $X(2)$ of
% $S_2$ is greater or equal to 3.
 Then there is an infinite operad $\p$ generated by $X$ such that every its element $x\in \p$
  is strongly nilpotent, that is, the suboperad generated by $x$ is finite.
\end{cor}

\rema In the language of varieties, the above Corollary looks as
follows:  there exist a variety of algebras with two binary
operations such that every operation in these algebras is
nilpotent, that is,
 for every operation (that is, homogeneous element of free algebra) there is a number $N$
 such that every composition of $N$ operations of this kind is zero for all
 possible substitutions of variables.

\begin{proof}[Idea of proof]
For every multilinear operation $p$ (say, $n$-ary, where $n \ge
2$) and sufficiently large number $d$, we define a suitable set of
relations $S(p,d)$ as a set  of all possible compositions of the
operation $p$ with itself of the following kind: in each of $d$
copies of $p$ in the compositions, at least $n-1$ of the inputs of
$p$ are replaced by variables. That is, every element of $S(p,d)$
looks as a ``branch'' of length $d$ whose nodes are market by $p$.
 Then we choose $d$ sufficiently large for each $p$ and use
Corollary~\ref{cor:GS_binary}.
\end{proof}

\end{document}